\documentclass[12pt,a4paper]{article}
\usepackage[centertags]{amsmath}
\usepackage{amsfonts}
\usepackage{amssymb}
\usepackage{amsthm}
\usepackage{newlfont}
\usepackage{indentfirst}
\title{An ILB- Manifold Structure on the Set of Riemannian
Metrics on a Noncompact Manifold}
\author{Catalin C. Vasii}
\date{}
\theoremstyle{plain}
\newtheorem{te}{Theorem}[section]

\newtheorem{pr}[te]{Proposition}
\newtheorem{co}[te]{Corollary}
\theoremstyle{definition}
\newtheorem{de}[te]{Definition}
\newtheorem{ex}[te]{Example}
\theoremstyle{remark}
\newtheorem{re}[te]{Remark}

\begin{document}
\maketitle
\typeout{:?0000}
 \begin{abstract}

 In this paper, using the structures of cone and bicone fields on
 vector bundles, the author introduces a ILB (inverse limit of
 Banach)- manifold structure on $\mathcal M$ the space of
 Riemannian metrics on a noncompact manifold $M$. In the last
 section, it is proven that, this way, on the open submanifold $\mathcal{M}_{finite}$
 of finite volume metrics, the canonical Riemannian metric is
 defined.

\par{\bf key words and frases:}  spaces of metrics, non compact manifolds;
\par {\bf msn subj. class:} 58D17.

 \end{abstract}
 \section{Preliminaries}
         First, let $M$  be a topological manifold, paracompact, with
  $\partial M= \emptyset$, which need not to be compact.
Let $(E,p,M)$ be a topological vector bundle over $M$.
\begin{de}\label{fvconuri} \cite{pa1}\\
A \textit{cone field} on the vector bundle $(E,p,M)$ is a map
$K:M\to \mathcal{P}(E),$\linebreak $x\mapsto K(x)\subset E_x$
which satisfies the following two conditions:\\ (K1)
$(\forall)x\in M,\ K(x)$ is a convex cone, closed in  $E_x$,
pointed, solid;\\ (K2) $\cup_{x\in M}Int(K(x))$\ and $\cup_{x\in
M}(E_x \backslash K(x))$ are open in $E$.
\end{de}
  In the definition  above,  a \textit{convex cone} is,
following \cite{Kras:ro}, a set $K$ which satisfies $K+K \subset
K$ and $(\forall)\lambda \geq 0,\ \lambda K \subset K.$ A cone $K$
which satisfies $K\cap -K= \{0\}$ will be called \textit{pointed
cone}, and \textit{a solid cone} is a cone which has interior
points in the topology of $E_x$.
\par The structure consisting by a vector bundle $(E,p,M)$   and a cone field  $K$ on it
is denoted by  $[(E,p,M);K]$.
\begin{ex}\label{metrici}\cite{pa2}
Let us consider now the bundle $(\mathcal{S}^2T^*M,p,M)$ of 2
times covariant symmetric  tensors on a given manifold $M$. We put
$(\forall) x\in M$
$$K_k(x):=\{t_x\in \mathcal{S}^2T^*M_x|r(t_x)=i_p(t_x)\},$$
where $r$ denotes the rank and $i_p$ denotes the positive inertia
index. Then $x\mapsto K(x):=\cup_{k=1}^n K_k(x)$ defines a cone
field on the bundle $(\mathcal{S}^2T^*M,p,M)$.
\end{ex}

\par There are local and global properties of  this structure,
exposed in \cite{pa1}.
\par  Consider now $\Gamma^0 (E)$, the space of continuous sections
of the bundle $(E,p,M)$.
 \begin{de}\label{secpoz} \cite{pa1}\\
 We call a \textit{positive section} of the structure
 $[(E,p,M);K]$, a section $\sigma \in \Gamma^0 (E)$ for which
 $\sigma (x)\in K(x),\ (\forall) x\in M$.
 \end{de}
 \par The set of positive sections is denoted by $K^0_\Gamma$;\  if
 on the space $\Gamma^0 (E)$ is considered
 the graph topology $WO^0$ then we have:
 \begin{pr}\label{codinsec}\cite{pa1} \\
 The set $K^0_\Gamma $ is a convex cone, pointed, solid, in
 $\Gamma^0 (E)$. Moreover, $K^0_\Gamma $ is a generating cone in
 $\Gamma^0 (E)$, i.e. $\forall \sigma \in \Gamma^0 (E)\  (\exists)\zeta_1,
\zeta_2 \in K^0_\Gamma$
 such that \linebreak $\sigma = \zeta_1 - \zeta_2$.
 \end{pr}
 The cone $K^0_\Gamma$ defines a partial order relation on $\Gamma^0 (E)$
 by
\begin{equation}\label{relordine}
 \sigma_1 , \sigma_2 \in \Gamma^0 (E),\ \sigma_1 \leq \sigma_2
 \stackrel{def}{\iff}\sigma_2 -\sigma_1 \in K^0_\Gamma.
\end{equation}
 \begin{pr}\verb"Papuc"\cite{pa1}\\
The pair $(\Gamma^0(E),\leq)$ is an ordered vector space, directed
on both sides. Endowed with the $WO^0$ topology is a topological
vector space iff $M$ is a compact manifold.
\end{pr}
\par Given a fixed $\zeta \in Int(K^0_\Gamma)$, we denote by $\Gamma^0_\zeta$
the set of  \linebreak $\zeta-$ measurable elements of $\Gamma^0
(E)$:
\begin{equation}\label{masurabil}
\Gamma^0_\zeta := \{\sigma \in \Gamma^0 (E)\mid (\exists)\lambda
\in \mathbb R_+ : -\lambda \zeta \leq \sigma \leq \lambda \zeta\},
\end{equation}
 and
we consider the map
\begin{equation}\label{normasect}
|\cdot |^0_\zeta :\Gamma^0_\zeta \to \mathbb R, \ |\sigma |_\zeta
:=min_{\lambda \in \mathbb R_+} \{-\lambda \zeta \leq \sigma \leq
\lambda \zeta\}.
\end{equation}

\begin{pr}\label{C} \cite{pa1}\\
1.\ $\Gamma^0_\zeta = \Gamma^0 (E)$ iff $M$ is a compact manifold;\\
2.\ The map $|\cdot |^0_\zeta$ defined by equation
(\ref{normasect}) from above is a norm on $\Gamma^0_\zeta$;\\
3.\ The set of all $\Gamma^0_\zeta$ is a covering of $\Gamma^0 (E)$;\\
4.\ If $\Gamma^0_c (E)$ denotes the subspace of compact support
sections, then \linebreak $\Gamma^0_c (E)\subset \Gamma^0_\zeta,
(\forall)\zeta \in Int(K^0_\Gamma) $;\\ 5.\ If $\zeta \in
\Gamma^0_{\zeta_1},$ with $\zeta, \zeta_1 \in Int (K^0_\Gamma)$
then $\Gamma^0_\zeta \subset \Gamma^0_{\zeta_1}$.
\end{pr}
\begin{te}\cite{vasii} \label{banach}
If $\zeta \in K^0_\Gamma$, then $(\Gamma^0_\zeta,
|\cdot|^0_\zeta)$ is a Banach space.
\end{te}
We consider now $(E,p,M)$, a $\mathcal C^k$- differentiable bundle
over $M$, a $\mathcal C^k$- differentiable manifold, $k\geq1$,
which need not to be compact.
\begin{de}\cite{vasii2}\label{supercone}
A {\it bicone field} on a vector bundle $(E,p,M)$ is the structure
consisting of a cone field $K$ on the bundle $(E,p,M)$ and a
second cone field $K_{TM}$ on the tangent bundle $(TM,p,M)$.
\end{de}
We will denote by $[(E,p,M);K;K_{TM}]$ the structure consisting of
a  bicone field on the vector bundle $(E,p,M)$.
\par
 The existence of a bicone field on a
vector bundle $(E,p,M)$ is equivalent with the existence of a non
zero section $\zeta \in \Gamma^0(E)$ and  of a nonzero vector
field on $M$.

\par Now, as a consequence of   the vector bundle
isomorphism
\begin{equation}\label{izopalais}
J^kE \cong \oplus_{i=1}^k\mathcal L_s(TM^i,E)
\end{equation}
from \cite{palais} page 90, we have
\begin{pr}\cite{vasii2}
If $[(E,p,M);K;K_{TM}]$ is a $\mathcal C^p$- differentiable vector
bundle endowed with a  bicone field then the vector bundle
$(J^kE,p,M)$ is endowed in a natural way with a cone field $K^k,\
(\forall)k\leq p$.
\end{pr}
Next, as usually, we will denote by convention $J^0E:=E,\
j^0\zeta:=\zeta$.
\begin{de}\cite{vasii2}\label{secposk}
A section $\zeta\ \in \Gamma^k(E)$ which satisfies $\zeta(x)\in
K(x)$ and $j^i\zeta(x)\in K^i(x),i=\overline{0,k}, (\forall)x\in
M$ will be called section {\it positive up to order $k$}.
\end{de}
The set of positive sections up to order $k$ will be denoted by
$K^k_\Gamma$. On $\Gamma^k(E)$, the space of $\mathcal C^k$-
differentiable sections we will consider the Whitney $WO^k$-
topology, which on a space of sections can be given by a base of
neighborhoods $W(\sigma_0, U)$, where $\sigma_0\in \Gamma^k(E)$
and $U$ is an open neighborhood of $Im(j^k\sigma_0)$ in $J^k(E)$:
\begin{equation}\label{bazak}
W(\sigma_0,U):=\{\sigma\in \Gamma^k(E)\mid j^k\sigma(x)\in
U,(\forall)x\in M\}.
\end{equation}
\begin{pr}\label{vasii2}
$K^k_\Gamma$ is a convex cone, closed, pointed and solid in the
space $(\Gamma^k(E),WO^k)$.
\end{pr}
\begin{co}\label{relordk}\cite{Kras:ro}
The cone $K^k_\Gamma$ defines on $\Gamma^k(E)$ an order relation
by $\sigma_1 \leq \sigma_2 \stackrel{def}{\iff}\sigma_2-\sigma_1
\in K^k_\Gamma$. In particular, this relation is directed on both
sides.
\end{co}
Let $\zeta \in Int(K^k_\Gamma)$ be fixed.
\begin{de}\cite{vasii2}\label{mask}
A section $\sigma \in \Gamma^k(E)$  for which exists $\lambda \in
\mathbb{R}_+$ s.t.
$$ - \lambda j^i \zeta \leq j^i \sigma \leq \lambda j^i \zeta,\ i=\overline{1,k}$$
will be called {\it $\zeta$- measurable up to order $k$}.
\end{de}
\par As in \cite{pa1}, we have that the map
$$|\cdot|^k_\zeta :\Gamma^k_\zeta \to \mathbb{R}_+,\
|\sigma|^k_\zeta:=\min\{\lambda \in \mathbb{R}_+\mid - \lambda j^i
\zeta \leq j^i \sigma \leq \lambda j^i \zeta,\
i=\overline{1,k}\}$$ is a norm on the vector space
$\Gamma^k_\zeta$ of $\zeta$- measurable sections up to order $k$,
and with this norm, $\Gamma^k_\zeta$ becomes a Banach space (the
proof is absolutely similar to the one from \cite{vasii}). The
open ball in the norm $|\cdot|^k_\zeta$, centered in $\sigma$, of
radius $\epsilon$, will be denoted by $B^k_\zeta(\sigma,\epsilon)$
and as in \cite{pa1}, coincides with the open centered intervals
in the order relation from corollary \ref{relordk}.
\par Let us denote now by $\tau^k$ the topology on $\Gamma^k(E)$
obtained by taking the path connected components of the $WO^k$-
topology.
\begin{te}\cite{vasii2}\label{topck}
For all $k\in \mathbb{N}$, the $\tau^k$- topology on $\Gamma^k(E)$
is the topology for which a basis of neighborhoods is given by
$$\{B^k_\zeta (\sigma, \epsilon)\mid \zeta \in Int(K^k_\Gamma),\
\sigma \in \Gamma^k_\zeta,\ \epsilon \gneq 0\}.$$
\end{te}

 \section{The ILB- manifold Structure on the Space of Riemannian Metrics}
   Let now $(E,p,M)$ be a smooth vector bundle, endowed with a bicone
field defined by the cone fields $K, K_{TM}$.
\begin{de}\cite{vasii2}\label{indefinit}
A smooth section $\zeta \in \Gamma(E)$ which satisfies
$\zeta(x)\in K_\Gamma ^k,\ (\forall) k\in \mathbb N $ will be
called a {\it indefinitely positive section}.
\end{de}
We will denote by $K_\Gamma$ the set of indefinitely positive
sections.
\begin{pr}\cite{vasii2}
The set $K_\Gamma$ is a  (nonempty) pointed closed convex cone in
$(\Gamma(E), WO^\infty)$.
\end{pr}
\begin{co}\cite{vasii2}
On $\Gamma(E)$ there is an order relation defined by $\sigma\leq
\sigma' \stackrel{def}{\iff} \sigma' -\sigma \in K_\Gamma(E)$.
\end{co}
Let $\zeta  \in \cap_k Int_{WO^k}K^k_\Gamma$ (this set is
nonempty, see \cite{vasii2}).
\begin{de}\cite{vasii2}
A section $\sigma \in \Gamma(E)$ which is $\zeta$- measurable
$(\forall)k\in \mathbb N$  will be called an {\it indefinitely
$\zeta$- measurable section}.
\end{de}
The set $\Gamma_\zeta (E)$ of indefinitely $\zeta$- measurable
sections is nonempty (e.g. $\zeta \in \Gamma_\zeta$ ) and is a
vector space.
\begin{pr}\cite{vasii2}
The space $\Gamma_\zeta (E)$ is the projective limit of the Banach
spaces $\Gamma^k_\zeta(E)$.
\end{pr}
\begin{co}\cite{vasii2} The following assumptions hold:\\
(i)\ $\Gamma_\zeta(E)$ is a complete, locally convex space;\\
(ii)\ The $\tau^\infty$- topology on $\Gamma(E)$ is the topology
for which a base of neighborhoods is given by the set
$$\{B^k_\zeta(\sigma, \epsilon)|\ \zeta  \in \cap_k Int_{WO^k}K^k_\Gamma,\ k\in
\mathbb{N},\ \epsilon \gneq 0 \};$$ (iii)\ The set
$\{\Gamma_\zeta(E),\Gamma^k_\zeta(E)\mid\ k\in N(0) \}$ is a ILB
(inverse limit of Banach)- chain, following Omori's definition
\cite{omori}, page 5.
\end{co}
\par Since in the infinite dimensional geometry the notion of
manifold might vary, we will refer in this paper to the notion
from \cite{mi-kri}, page170, for which the differences from the
finite dimensional correspondent is that for each chart is allow a
different model space, and the chart changing is require  to be
only smooth instead of smooth diffeomorphism.
\begin{te}$\Gamma (E)$ is a smooth manifold modelled by the
ILB-spaces $\Gamma_\zeta (E)$.
\end{te}
\proof From \cite{vasii} and \cite{vasii2} we have
$\Gamma(E)=\underrightarrow{\lim_{\zeta}}\underleftarrow{\lim_k}\Gamma^k_\zeta$.
The topology induced above  on $\Gamma(E)$   is the $\tau^\infty$-
topology. Then, again by the equation above,
$\Gamma(E)=\cup_{\zeta \in Int(K_\Gamma)}\Gamma_\zeta(E)$.
\par Let $\sigma_0\in  \Gamma(E)$. There exists  a positive section $\zeta_0 \in
Int(K_\Gamma)$ such that $\sigma_0\in \Gamma_{\zeta_0}(E)=\cap_k
\Gamma^k_{\zeta_0}(E)$. Obviously, $U_{\zeta_0}(\sigma_0):=\cap_k
B^k_{\zeta_0}(\sigma_0)$ is a nonempty open in $\tau^\infty$-
topology neighborhood of $\sigma_0$. Let
$\phi_{\sigma_0}:U_{\zeta_0}(\sigma_0)\subset \Gamma(E)\to
\Gamma_{\zeta_0}$ be the restriction of the identity map
$Id_{\Gamma_{\zeta_0}(E)}$. The pair $(U_{\zeta_0}(\sigma_0),
\phi_{\zeta_0})$ is a chart around $\sigma_0$.
\par The charts changing is smooth. Indeed, Let
$(U_{\zeta_1}(\sigma_1), \phi_{\sigma_1}),(U_{\zeta_2}(\sigma_2),
\phi_{\sigma_2})$ be two charts with $U_{\zeta_1}(\sigma_1)\cup
U_{\zeta_2}(\sigma_2)\neq \emptyset$. In particular, it follows
that $U_{\zeta_1}(\sigma_1)\cup U_{\zeta_2}(\sigma_2)\subset
\Gamma_{\zeta_1}\cap \Gamma_{\zeta_2}$. But from \cite{pa1}, the
set $\{\Gamma_\zeta(E) | \zeta \in Int(K_\Gamma)\}$ is ordered and
directed on both sides, by the inclusion. So there exists
$\zeta_0\in Int(K_\Gamma)$ such that $U_{\zeta_1}(\sigma_1)\cap
U_{\zeta_2}(\sigma_2)\subset \Gamma_{\zeta_1}\cap
\Gamma_{\zeta_2}\subset \Gamma_{\zeta_0}(E)$, and so the chart
changing $\phi_{\sigma_2} \circ \phi^{-1}_{\sigma_1}$ is the
restriction to an open set of the identity map
$Id_{\Gamma_{\zeta_0}(E)}$, and so is smooth. \hfill Q.E.D.\\
 \begin{re}\label{vartens}
 In virtue of the example \ref{metrici}, $\Gamma
 (\mathcal{S}^2T^*M)$, the space of two times covariant, symmetric
 tensor fields on the manifold $M$ has the structure of a
 ILB-manifold, modelled by the spaces
 $\Gamma_g(\mathcal{S}^2T^*M)$, with $g\in Int(K_\Gamma)$.
 \end{re}
From \cite{michor} $\mathcal{M}=Int(K_\Gamma)\cap \Gamma
 (\mathcal{S}^2T^*M)$, the space of all Riemannian metrics on
 the manifold $M$ is $\tau^\infty$- open in $\Gamma
 (\mathcal{S}^2T^*M)$.
\begin{co}\label{varmetr}
The space $\mathcal M$ of all Riemannian metrics on $M$ is an open
submanifold of $\Gamma (\mathcal{S}^2T^*M)$.
\end{co}

 \section{The Riemannian Geometry of the Space of Riemannian Metrics of Finite Volume}
   We denote by $\mathcal M_{finite}$ the set of  all Riemannian
metrics of finite volume on $M$.
\begin{re} $\mathcal{M}_{finite}$ is $\tau^\infty$- open in $\mathcal
M$. Indeed, let $(g_n)_{n\geq 0}$ a sequence of Riemannian metrics
that converges in the $\tau^\infty$- topology to $g_0$, a finite
volume metric. In particular, it follows that $(\forall) n\geq 0,
g_n $ and $g_0$ differ only on a compact set, so each $g_n$ is a
finite volume metric.
\end{re}
\par On $\mathcal M$ there is a canonical Riemannian metric $G$,
invariant under the natural action by pull- back of  the group
$Diff(M)$ of diffeomorphisms of $M$  on $\mathcal M$, described in
\cite{ebin}, or \cite{michor}:
\begin{equation}
G_g:T_g\mathcal M \times T_g \mathcal M\to \mathbb R,\
G_g(h,k)=\int_M trace (g^{-1}hg^{-1}k)d\nu_g,
\end{equation}
To make clear the notation  $g^{-1}hg^{-1}k$ we can regarde the
bundle $\mathcal{S}^2 (T^*M)$ as $\{h\in \mathcal{L}(TM,T^*M)|\
h^t=h\}$, subbundle of $\mathcal{L}(TM,T^*M)$, where $h^t$ is the
composition $TM\stackrel{i}{\hookrightarrow}T^{**}M
\stackrel{h^*}{\to}T^*M$. On the other side, since $g\in \mathcal
M$, as a Riemannian metric is a fiberwise inner product on $TM$ it
induces a fiberwise inner product on any tensor bundle over $M$,
in particular on $\mathcal{S}^2 (T^*M)$. This is, in fact
$<\cdot,\cdot>=trace(g^{-1}\cdot g^{-1}\cdot)$. For the metric
$G_g$, instead of the notation above, we will use the classical
notations from Riemannian geometry (the '$\sharp$' symbol demotes
the 'sharp' isomorphism induced by the metric $g$ so we will omite
to put indices as $\sharp_g$ ):
$$ G_g:T_g\mathcal M \times T_g \mathcal M\to \mathbb R,\
G_g(h,k)=\int_M \sum_{i=1}^n h(k(E_i)^\sharp, E_i)d\nu_g,$$ Where
$(E_i)$ denotes a local field of orthonormal frames.
\begin{te}
The Riemannian metric $G_g$ is defined on the tangent space $T_g
\mathcal{M}_{finite}=\Gamma_g$.
\end{te}
\proof Since $h\in \Gamma_g$, we have $h\in \Gamma^0_g$. This
means that $(\exists )\lambda \in \mathbb{R}_+$ s.t. $\lambda
g\leq h\leq \lambda g $. Because of equation (\ref{normasect}), we
have that $-|h|_g^0 g\leq h\leq |h|_g^0 g$. Hence, as in
\cite{pa1}, in particular,
$$-|h|_g^0 g(k(E_i)^\sharp ,E_i)\leq h(k(E_i)^\sharp ,E_i)\leq |h|_g^0
g(k(E_i)^\sharp ,E_i),\ i=\overline{1,n};$$ By summation, we have
$$-|h|_g^0\sum_{i=1}^n g(k(E_i)^\sharp ,E_i)\leq \sum_{i=1}^n h(k(E_i)^\sharp ,E_i)
\leq |h|_g^0 \sum_{i=1}^ng(k(E_i)^\sharp ,E_i),$$ and this means
\begin{equation}\label{unu}
-|h|_g^0 \sum_{i=1}^n k(E_i,E_i)\leq \sum_{i=1}^n h(k(E_i)^\sharp
,E_i)\leq |h|^0_g \sum_{i=1}^n k(E_i,E_i).
\end{equation}
But $k\in \Gamma_g$, so we have $k\in \Gamma^0_g$. This means that
$(\exists )\lambda \in \mathbb{R}_+$ s.t. $\lambda g\leq k\leq
\lambda g $. As above, $-|k|_g^0 g\leq k\leq |k|_g^0 g$, and in
particular
\begin{equation}\label{doi}
-|k|_g^0 g(E_i, E_i) \leq k(E_i,E_i)\leq |k|_g^0 g(E_i,E_i),\
i=\overline{1,n};$$ By summation
$$-|k|_g^0 \sum_{i=1}^n g(E_i, E_i) \leq \sum_{i=1}^n k(E_i,E_i)\leq |k|_g^0
\sum_{i=1}^n g(E_i,E_i).
\end{equation}
From equations (\ref{unu}) and (\ref{doi}) follows that
$$-|h|_g^0 |k|_g^0 \sum_{i=1}^n g(E_i,E_i)\leq \sum_{i=1}^n h(k(E_i)^\sharp
,E_i)\leq |h|^0_g |k|_g^0 \sum_{i=1}^n g(E_i,E_i)$$ and so
$$-n|h|_g^0 |k|_g^0 \leq \sum_{i=1}^n h(k(E_i)^\sharp
,E_i)\leq n|h|^0_g |k|_g^0. $$ Now, by integrating with respect to
the measure $\nu_g$
$$-n|h|_g^0 |k|_g^0 Vol(M,g)\leq G_g(h,k)\leq n|h|_g^0 |k|_g^0 Vol(M,g)$$
\hfill Q.E.D.

 \begin{tabbing}
 C\u{a}t\u{a}lin C. Vasii\\
 Departament of Mathematics and Computer Sciences, \\
 The West University of  Timisoara,\\
 Bv. V. P\^arvan nr.4, 300.223 Timisoara, Rom\^ania,\\
 catalin@math.uvt.ro\\
 \end{tabbing}
\end{document}